\documentclass[12pt]{article}

 \usepackage{
 amsmath,amsfonts,
 amssymb,amsthm
}

 \numberwithin{equation}{section}
 \newtheorem{ex}{Example}[section]
 \newtheorem{prop}{Proposition}[section]

 \newtheorem{cor}{Corollary}[section]
\newtheorem{lem}{Lemma}[section]
 \newtheorem{dfn}{Definition}[section]

 \newtheorem{rmkk}{Remark}[section]

 \newcommand{\p}{\partial}

 \newcommand{\sG}{\mathcal G} 

 \newcommand{\sF}{\mathcal F}

 \newcommand{\al}{\alpha}
 \newcommand{\be}{\beta}
 \newcommand{\de}{\delta}

 \newcommand{\ga}{\gamma}

 \newcommand{\De}{\Delta}
 \newcommand{\Ga}{\Gamma}
 \newcommand{\La}{\Lambda}
 \newcommand{\Om}{\Omega}
 \newcommand{\la}{\lambda}
 \newcommand{\Si}{\Sigma}
 \newcommand{\Th}{\Theta}

 \newcommand{\R}{\mathbb R}
 \newcommand{\Z}{\mathbb Z}

 \newcommand{\rk}{\operatorname{rank}}

 \renewcommand{\Im}{\operatorname{Im}}

\newcommand{\Ker}{\operatorname{Ker}}

\begin{document}

\title{Geometry of unimodular systems}
  \author{I.V.Artamkin \thanks{Faculty of Mathematics, National Research University Higher School
of Economics, Usacheva str., 6, 119048 Moscow, Russian Federation
\endgraf
MIREA - Russian Technological University, Vernadskogo Pr., 78,
Moscow 119454, Russia}}
 \date{}
 \maketitle
 \begin{abstract}  A  collection of vectors in a real vector space is called
a {\it unimodular system} if  any of its maximal linearly
independent subsets   generates the same free abelian group. This
notion is closely connected with {\it totally unimodular matrices}:
rows or columns of a totally unimodular matrix form a unimodular
system and the matrix of coefficients of  expansions of all vectors
of a unimodular system with respect to its  maximal linearly
independent subset  is totally unimodular.

In this paper we show that a unimodular system defines the following
geometric data: a Euclidean space, an integral lattice in it, and a
reflexive lattice zonotope. The discriminant of the lattice is equal
to  the number of maximal linearly independent subsystems, and we
call this number the complexity of the unimodular system.  For a
unimodular system   $\Omega $  we also define the  Gale dual
unimodular system $\Omega ^{\bot}$ which has the same complexity.
These notions may be illustrated by the well-known graphic and
cographic unimodular systems of a graph. Both graphic and cographic
unimodular systems have the same complexity which is equal  to the
complexity of the graph. For graphs without loops and bridges the
graphic and the cographic unimodular systems  are Gale dual to each
other.

We describe this geometric data for certain examples: for the
graphic and the cographic unimodular systems of a generalized
theta-graph, consisting of two vertices connected by $N$ edges, for
the cographic system of the complete graph $K_N$,  and 
for  the
famous Bixby-Seymour
 unimodular  system, which is neither graphic nor cographic.
\end{abstract}

 \section{Introduction}

A real matrix is called totally unimodular  (see \cite{Schr} for
details) if all of its minors take one of the values 0, 1 or -1.
This implies that every maximal linearly independent set of its rows
(or columns) generates over $\Z $ the same free abelian group. This
observation suggests a more geometric approach to totally unimodular
matrices: if any maximal linearly independent subset of some
collection of vectors generates the same free abelian group, then
the matrix of coefficients of  expansions of all vectors of this
system with respect to this basis is totally unimodular (see
Proposition \ref{prop_sys-matr} below).  We will call such
collections of vectors {\it unimodular systems} of vectors.   We
shall consider vectors in unimodular systems up to a sign; one
vector may appear
  several times in a unimodular system.
For our purposes it is more convenient to pass to the dual space and
consider unimodular systems of linear forms --- see Definition
\ref{def_unumod_sys} below.

  In this paper
  we show that an
 unimodular system $\Om $ defines the following geometric data: a
 real Euclidean vector
  space $W_{\Omega }$,   an integral  lattice $L_{\Omega } \subset W_{\Omega }$
dual to the one defined by the linear forms from the unimodular
system $\Om $, and a polytope  $ \De _{\Omega }\subset W_{\Omega }$,
whose facets are defined by the linear forms from $\Om $.  We prove
that the discriminant of the lattice $L_{\Om }$ is equal to the {\it
complexity}   $c(\Om )$ of the unimodular system $\Omega $
 (see Definition \ref{def_complexity}) which
is the number of maximal linearly independent subsystems (i.e.
bases)  in $\Om $ (see Proposition \ref{Kirh} below). In Section
\ref{section_polytope} we prove that $\De  _{ \Om  } $ is a
reflexive lattice polytope and also a zonotope. The polytope $\De _{
\Om } $ first appeared in an algebro-geometric context in \cite{A1}:
the generalized Jacobian of a maximally degenerated stable curve is
a toric variety defined by the polytope corresponding to the graphic
unimodular system  of the dual graph of the curve.

In  Section \ref{section_duality} for an unimodular system $\Om $ we
define the  {\it  Gale dual} unimodular system $\Omega ^{\bot}$ of
the same complexity. The double
   dual system $(\Om  ^{\bot}) ^{\bot}$ is a subsystem of $\Om $ and
   coincides with $\Om $ if and only if  $L_{\Om  }$ has no vectors
   of length 1.

The above notions may be illustrated by the well-known graphic and
cographic unimodular systems of a graph. We recall the construction
of these unimodular systems in Section  \ref{section_graph}. The
complexity of both systems is equal  to the complexity of the graph
and these systems are Gale dual to each other for  graphs without
loops and bridges.

We describe the lattice  $L_{\Om  }$ and the polytope $\Delta _{ \Om
}$ for the following particular cases.
\begin{itemize}
  \item For a generalized
theta-graph 
(two vertices connected by $N$ edges) the lattice for the
cographic system is the one-dimensional lattice $\sqrt N \ \Z $
and the polytope is the segment $[-\sqrt N ; \sqrt N ]$. See
Example \ref{ex_poly_gen_theta_graph1}.
  \item The lattice for the graphic system for a generalized
theta-graph is the root lattice  of type $A_{N-1}$, namely $ \{
 (x_1,\ldots , x_N)\in \Z ^N \ , \ \sum _{i=1}^N x_i =0 \} $,
and the polytope is the
 intersection of two simplices that are centrally symmetrical to
each other, one of which is given by its vertices $P_1, \ldots ,
P_N$ where $P_k=(-1,-1, \ldots, -1,  N-1, -1, \ldots, -1)$ (all
the components are $-1$, except for the $k$-th component, which
is $N-1$). See Example \ref{ex_poly_gen_theta_graph2}.
  \item The lattice for the cographic system  for the
 complete graph    $K_N$ on $N$ vertices is $ N\cdot A_{N-1}^*$,
 where  $A_{N-1}^*$ is the lattice dual to the root lattice of
 type  $A_{N-1}$, see Corollary
 \ref{cor_A_N_*_for_complete_graph}. The polytope is the
 projection of the standard cube $[0;1]^N$ parallel to its main
 diagonal, see Example \ref{poly_K_N}. In particular for $N=4$
it is  the rhombic dodecahedron.
\end{itemize}

Finally in Proposition \ref{prop_Seymour} we describe the geometry
of the famous Bixby-Seymour
 unimodular  system
(\cite{Bixby}, \cite{Seym}) which is neither graphic nor cographic.
The Bixby-Seymour
 system is self-dual and its complexity is 162. The
corresponding polytope is the convex hull of the union of two
  regular simplices in Euclidean space of dimension 5 that are
  centrally symmetrical to each other.
The square of the position vector of each vertex is 10, and the
lattice is generated by the position vectors of  the midpoints of
the edges of these simplices, whose squares are 4. This implies that
the automorphism group of the Bixby-Seymour
  system is $\Z _2 \times S_6$.

\section{Totally unimodular matrices and unimodular systems.}
Consider a totally unimodular matrix  $A$ of full rank having $n$
columns and $N$ rows where $N\ge n$. We shall consider rows of $A$
as coefficients of linear forms defined on some $n$-dimensional real
linear space $U$.
 Rows of any nonzero  $n\times n$-minor generate over $\Z
$ a free abelian subgroup of $U^*$ of rank $n$, and since all such
minors are equal $1$ or $-1$, this abelian group is exactly the same
for all such minors. The following definition is based on this
property.
\begin{dfn}\label{def_unumod_sys}
Fix a real $n$-dimensional  vector space $U$ and let $\Xi $ be a
collection\footnote{Since    $\Xi $ may have repeating elements, it
is, strictly speaking, not a set, but a multiset, so we use square
brackets rather then curly brackets to denote it.} of nonzero linear
forms on $U$: $\Xi =[ \xi _1, \ldots , \xi _N]$, $\xi _i \in U^*$.
We shall call  $\Om = (U,\Xi )$ {\bf unimodular system} if any
maximal linearly independent subset of $\Xi $ generates over $\Z$
the same free abelian group $M\subset U^*$ of rank $n$.
\end{dfn}
Note that the forms $\xi _i  $  are not necessarily all different,
and we will consider these forms up to  sign change.
\begin{dfn}\label{def_isom_unumod_sys} We shall call two unimodular
systems  $(U,\Xi )$ and   $(U',\Xi ')$ ($\Xi '=[ \xi _1', \ldots ,
\xi _{N'}']$) isomorphic  if $N=N'$ and there is a linear
isomorphism $f: U\to U'$ such that after appropriate  renumbering of
forms $f^* (\xi _i')=\pm \xi _i$.
\end{dfn}
\begin{dfn}\label{def_complexity}
  {\bf Complexity}   $c(\Omega )$ of   unimodular system $ \Om =(U,\Xi )$
 is the number of
 maximal  rank subsets of $\Xi $.
\end{dfn}

The following Proposition describes the above mentioned connection
between unimodular systems and totally unimodular matrices.

\begin{prop}\label{prop_sys-matr}
 a) Let $A=(a_{p,q})$ be a  $N\times n$ totally
unimodular matrix of maximal rank and $N\ge n$. For each row of $A$
consider linear function ${\al _k (x)=\sum _{i=1}^n a_{k,i}x_i}$
($x=(x_1, \ldots ,x_n)\in \R ^n$), $1\le k\le N$. Then $(\R ^n, [\al
_1,
\ldots ,\al _N])$ is an unimodular system.\\
b) Let  $ \Om = (U,\Xi )$ be an unimodular system, $\dim U = n$,
$\Xi =[ \xi _1, \ldots , \xi _N]$, $\xi _i \in U^*$, and let
$M\subset U^*$ be  the  free abelian group  generated by $\Xi $ over
$\Z$. Choose any basis $\xi _{i_1}, \ldots , \xi _{i_n} $ of $M$ and
expand $\xi _k=\sum _{s=1}^n a_{k,s}\xi _{i_s}$ for  $1\le k\le N$.
Then the matrix
\begin{equation}\label{A_Om}
A_{\Om }= (a_{p,q})
\end{equation}
 is totally unimodular.
\end{prop}
Following \cite{Schr} any maximal rank submatrix $C$ of $A$ will be
called  {\it basis} of $A$, so there is one-to-one correspondence
between the bases of matrix $A$ and those bases of vector space
spanned by rows of $A$, which consist of the rows of the  matrix
$A$. The rows of any basis $C$ of the matrix $A$ form a basis of the
vector space generated by rows of $A$ and the rows of matrix
$AC^{-1}$ present the expansion of the rows of $A$ over this basis.
Assumption $\det C= \pm 1$ in a) means that $C^{-1}$ is integer and
therefore all the rows of $A$ lie in the abelian subgroup generated
by the rows of the basis $C$, which proves a).

Note that we can renumber elements of  $\Xi $ in b) so that the
basis will be  $\xi _1, \ldots , \xi _n $ and then the matrix
$A_{\Om }$ will be
\begin{equation}\label{A_Om}
A_{\Om }=\left(
           \begin{array}{c}
             E \\
             \widetilde{A} \\
           \end{array}
         \right)
\end{equation}
where $E$ is the unit $n\times n$ matrix. Any base of the matrix
$A_{\Om }$ corresponds to base of $M$ and therefore the expansion of
any row over this basis is integer. Thus the matrix $A_{\Om } C^{-1}
$ is integer for any base $C$ of the matrix $A_{\Om }$. Theorem 19.5
from \cite{Schr} states that then the matrix $A_{\Om } C^{-1} $ is
totally unimodular for any base $C$ of  $A_{\Om }$. But for the base
consisting of $n$ first rows of  $A_{\Om } \ $ $C=E$ (see
(\ref{A_Om})), so $A_{\Om }$ is totally unimodular. $\square $

Now let us give the simplest  examples.
\begin{ex}\label{ex_Upsilon}
The most trivial unimodular system consists of one linear form $\xi
\in U^*$ for one-dimensional vector space $U$. Let us denote this
system by $\Upsilon $. $A_{\Upsilon }$ is then $1\times 1$ matrix
$(1)$. Its complexity $ c(\Upsilon )= 1$.
\end{ex}
\begin{ex}\label{ex_1N}
For $\dim U =1$ and any $N>1$ there is exactly one unimodular system
$(U, \Xi )$ with $|\Xi |=N$, namely the system $\Si _N =(U,
[\underbrace{\xi, \xi , \ldots , \xi }_{N \ \mbox {\scriptsize
times}}])$, it corresponds to $N\times 1$ matrix
\begin{equation}\label{1N}
    \left(
      \begin{array}{c}
        1 \\
        \vdots \\
        1 \\
      \end{array}
    \right) ,
\end{equation}
 Its complexity $ c(\Si _N )= N$.
\end{ex}
\begin{ex}\label{ex_n=2}
For $n=2$ there are exactly two unimodular systems without multiple
elements. They are defined  by  totally unimodular matrices
\begin{equation}\label{n=2_N=2}
     \left(
       \begin{array}{cc}
         1 & 0 \\
         0 & 1 \\
       \end{array}
     \right)
     , \quad N=2,
\end{equation}
and
\begin{equation}\label{n=2_N=3}
    \left(
      \begin{array}{cc}
        1 & 0 \\
        0 & 1 \\
        1 & 1 \\
      \end{array}
    \right)
      , \quad  N=3.
\end{equation}
For (\ref{n=2_N=2}) complexity is equal to   $1$, and for
(\ref{n=2_N=3}) it is $3$. Any unimodular system with multiple
elements and $n=2$ may be obtained from
 (\ref{n=2_N=2}) or  (\ref{n=2_N=3}) by repeating several times some
 of their rows.
\end{ex}
Basic classical  examples of unimodular systems come from graphs.
Each graph defines two unimodular systems named {\it graphic} and
{\it cographic}. We shall describe their constructions in the next
section.

Finally let's present the famous Bixby-Seymour unimodular  system
(\cite{Bixby}, \cite{Seym}) which is neither graphic nor cographic.
\begin{ex}\label{ex_Seymour}
 We start with matrix
\begin{equation}\label{preSeymour}
Q= \left(
      \begin{array}{ccccc}
        1 & 1 & 0 & 0 & 0 \\
        0 & 1 & 1 & 0 & 0 \\
        0 & 0 & 1 & 1 & 0 \\
        0 & 0 & 0 & 1 & 1 \\
        1 & 0 & 0 & 0 & 1 \\
        1 & 0 & 1 & 0 & 0 \\
        0 & 1 & 0 & 1 & 0 \\
        0 & 0 & 1 & 0 & 1 \\
        1 & 0 & 0 & 1 & 0 \\
        0 & 1 & 0 & 0 & 1 \\
      \end{array}
    \right)
\end{equation}
which presents all $\binom 52 =10$ ways to arrange 2 ones and 3
zeros in a row. $Q$ is not totally unimodular for all its maximal
minors are equal to $\pm 2$ or zero. But this fact means that any 5
linearly independent rows of $Q$ generate the same rank 5 free
abelian group and therefore we shall get a totally unimodular matrix
expanding all the rows by this basis. Let us choose the first 5 rows
for basis and then the matrix
\begin{equation}\label{Seymour}
   R=  Q\cdot \left(
      \begin{array}{ccccc}
        1 & 1 & 0 & 0 & 0 \\
        0 & 1 & 1 & 0 & 0 \\
        0 & 0 & 1 & 1 & 0 \\
        0 & 0 & 0 & 1 & 1 \\
        1 & 0 & 0 & 0 & 1 \\
        \end{array}
    \right)
 ^{-1}   =
\left(
      \begin{array}{ccccc}
        1 & 0 & 0 & 0 & 0 \\
        0 & 1 & 0 & 0 & 0 \\
        0 & 0 & 1 & 0 & 0 \\
        0 & 0 & 0 & 1 & 0 \\
        0 & 0 & 0 & 0 & 1 \\
        0 & 0 & 1 & -1 & 1 \\
        1 & 0 & 0 & 1 & -1 \\
        -1 & 1 & 0 & 0 & 1 \\
        1 & -1 & 1 & 0 & 0 \\
        0 & 1 & -1 & 1 & 0 \\
      \end{array}
    \right)
\end{equation}
is totally unimodular. We shall prove that it is  neither graphic
nor cographic and discuss its geometry
 in Proposition \ref{prop_Seymour}.
\end{ex}

\section{Graphic and cographic  unimodular
systems.}\label{section_graph}
Let $\Ga $ be a connected graph, $V (\Ga )$ is the set of vertices
of $\Ga $,  $E (\Ga )$ is the set of edges, $|E (\Ga )|=N$, $|V (\Ga
)|= m$. Recall standard definition of homology of graph $\Gamma $
considered as
 $1$-dimensional  simplicial  complex. First we fix an orientation on
the edges of $\Gamma $ and consider vector space of 0-chains
\begin{equation}\label{0-chains}
  C_0 (\Ga ,\R ) =\{ \sum _{  v \in V( \Gamma )}
a_v \cdot v, \quad a_v\in \R \}
\end{equation}
and vector space of 1-chains
\begin{equation}\label{1-chains}
 C_1 (\Ga ,\R ) =\{ \sum _{  e \in E( \Gamma )}
 b_e \cdot e, \quad b_e\in \R \} .
\end{equation}
$\dim C_0 (\Ga ,\R ) =m$, $\dim C_1  (\Ga ,\R )=N$.  The border map
\begin{equation}\label{border_map}
\p _1  (\Ga  ) : C_1 (\Ga ,\R ) \to C_0 (\Ga ,\R )
\end{equation}
is a linear map defined by its action  on the elements of the basis
of $C_1 (\R )$ by
\begin{equation}\label{diff}
 \p  _1  (\Ga  )(e)= \sum _{  v \in V( \Gamma )} i(e,v) \cdot v,
\end{equation}
where
\begin{equation}\label{i_e_v}
i(e,v) =  \left\{
\begin{array}{rl}
                   1 &  \mbox{ if under the chosen orientation the edge }  e
                   \mbox{ points to }  v  \\
                   -1 &   \mbox{ if under the chosen orientation the edge }  e
                   \mbox{ points from }  v  \\
                   0 &    \mbox{ if } e \mbox{ and  } v  \mbox{ are not incident}
                 \end{array}
   \right.
\end{equation}

Then  $\Ker \p   _1  (\Ga  )$ is just the  homology space  $H_1 (\Ga
 , \R)$. Note that if \linebreak
  ${\sum _{  e \in E( \Gamma )}
 b_e \cdot e \in \Ker \p   _1  (\Ga  )}$  then for every $v\in V(\Ga
 )$ holds the following relation
\begin{equation}\label{Kirh_eq}
\sum _{  e \in E( \Gamma )} i(e,v) b_e =0.
\end{equation}
This means that we may interpret the elements of   $\Ker \p   _1
(\Ga  )$ as flows on the graph  $\Ga  $, and the relation
(\ref{Kirh_eq}) is exactly the  Kirchhoff's condition for the vertex
$v\in V(\Ga )$.

We shall denote the space  $\Ker \p   _1 (\Ga  )$ by $F_{\Ga }$. So
 $F_{\Ga }$ is a linear subspace in  $C_1  (\Ga ,\R )$
 defined   by any $m-1$ of  Kirchhoff's  conditions (\ref{Kirh_eq}),
 which are linearly independent while the last condition is the sum of the
 the others.

Let us fix standard scalar products on $C_0  (\Ga ,\R )$ and $C_1
(\Ga ,\R )$ so that vertices in $C_0  (\Ga ,\R )$ and edges in $C_1
(\Ga ,\R )$ become orthonormal bases. Then $C_0  (\Ga ,\R )$ and
$C_1 (\Ga ,\R )$ may be identified with the dual spaces of cochains
${C^0 (\Ga ,\R )=(C_0 (\Ga ,\R ))^*}$ and $C^1 (\Ga ,\R )=(C_1 (\Ga
,\R ))^*$ and $\p   _1  (\Ga  )^*$ becomes the coborder operator
\begin{equation}\label{coborder_map}
\p _1  (\Ga  )^* : C_0 (\Ga ,\R ) \to C_1 (\Ga ,\R )
\end{equation}
which acts on the elements of the basis of $C_0 (\Ga ,\R )$ by
\begin{equation}\label{codiff}
 \p  _1  (\Ga  )^*(v)=  \sum _{  e \in E( \Gamma )} i(e,v)\cdot e .
\end{equation}
Thus $C_1  (\Ga ,\R )$ is the orthogonal direct sum
\begin{equation}\label{C_1=Oplus}
  C_1  (\Ga ,\R )=\Ker \p  _1  (\Ga  )\oplus \Im \p  _1  (\Ga  )^*.
\end{equation}

The matrix of the linear mapping $\p  _1  (\Ga  )^*$ with respect to
the chosen bases  is the matrix $I_{\Ga } $ defined by
\begin{equation}\label{incidence_matrix}
     (I_{\Ga }) _{p,q}=i(e_p, v_q).
\end{equation}
Note that $I$ is
 the matrix
 matrix  of incidence of the graph $\Ga $ with the chosen
orientation.

The mapping $  \p  _1  (\Ga  )^*$ in terms of cochains may be
interpreted as a gradient of discrete function defined on the
vertices of the graph. We shall denote the space $\Im \p  _1  (\Ga
)^*$ by $G_{\Ga }$, then the orthogonal decomposition
(\ref{C_1=Oplus}) takes the following form:
\begin{equation}\label{C_1=OplusFG}
  C_1  (\Ga ,\R )=F_{\Ga }\oplus G_{\Ga  }.
\end{equation}

While elements of $F_{\Ga } $ are linear combinations of cycles, the
 elements of $ G_{\Ga  }$ also have a similar geometric interpretation:
 they are  linear combinations of cuts. Recall that a cut is a 1-chain
which is the sum
 of edges pointing from $V'$ to $V''$ for a certain partition
 $V=V'\sqcup V''$; in other words it
 is $\p  _1  (\Ga  ) ^* (\sum  _{v \in V''}v)$.

Note that $\dim F_{\Ga  } =N-m+1$,  $\dim G_{\Ga  } = m-1$.

All the listed constructions do not depend on the initial choice of
orientation on the graph; changing the orientation of an edge only
changes the sign of the corresponding row of matrix $I_{\Ga }$.

We shall be interested in integer points in the described spaces.
Consider the standard lattice in $C_1(\Ga  ,\R )$:
\begin{equation}\label{C_1(Z)}
C_1(\Ga  ,\Z ) =\{ \sum _{  e \in E( \Ga ) }
 b_e \cdot e, \quad b_e\in \Z \} \subset C_1(\Ga  ,\R ).
\end{equation}
\begin{prop}\label{lattices_for_graph}
$F_{\Ga } \cap C_1(\Ga  ,\Z )$ and  $G_{\Ga } \cap C_1(\Ga ,\Z )$
are lattices of maximal rank in $F_{\Ga }   $ and in $G_{\Ga } $.
\end{prop}
Fix a spanning tree $T$ in $\Ga $.

Then for any edge $e$ such that $e\notin T$, there is a unique path
in $T$ connecting the vertices of $e$; adding $e$ we get a cycle
$\ga _e \in F_{\Ga } \cap C_1(\Ga  ,\Z )$. (We choose the same
orientation on  $\ga _e$ as was initially fixed on $e$.) Now it is
not hard to see that
\begin{equation}\label{basis_ga}
 \{ \ga _e \ , \ e\notin T \}
\end{equation}
is a basis of $ F_{\Ga } $ and of $ F_{\Ga }  \cap C_1(\Ga  ,\Z )$:
if $\be = \sum _{e\in E} b_e e \in \Ker \p $ then
 $\be = \sum _{e\notin T} b_e \ga _e$, because
 the difference $ (\sum _{e\in E} b_e e) -
 ( \sum _{e\notin T} b_e \ga _e )$ is a cycle on
 the tree $T$ and therefore it is zero.

Similarly we construct a basis of  $  G_{\Ga } $ and of  $  G_{\Ga }
\cap C_1(\Ga  ,\Z )$  using edges of $T$. Any edge $e\in T$
  disconnects the vertices of $T$ into $V=V' \sqcup V''$, where $e$
  under fixed orientation points from $ V' $ to $ V''$. Then the cut $\de _e
  =\p ^* (\sum _{v \in V''} v) \in   G_{\Ga }  \cap C_1(\Ga  ,\Z )$ and
\begin{equation}\label{basis_de}
\{   \de _e \ , e\in T \}
\end{equation}
is the  bases of  $    G_{\Ga }  $ and of  $ G_{\Ga }  \cap C_1(\Ga
,\Z )$: if $\be = \sum _{e\in E} b_e e \in \Im \p ^* $ then
 $\be = \sum _{e\in T} b_e \de _e$. $\square $

We shall denote the lattices $F_{\Ga }\cap C_1(\Ga  ,\Z )$ and
 $G_{\Ga }\cap C_1(\Ga  ,\Z )$ respectively by $F_{\Ga } ( \Z )$ and
 $G_{\Ga } ( \Z )$.

Using the scalar product we may consider edges of $\Gamma $ as
linear functions on $C_1( \Ga ,\R )$ and therefore on  $F_{\Ga } $
and $G_{\Ga } $. An edge $e$ is a zero function on $F_{\Ga }  $ if
and only if it is a bridge (i.e. removing $e$ makes graph $\Ga $
disconnected).
 An edge $e$ is a zero function on $G_{\Ga } $ if and only if it
is a loop. Denote by $E_b\subset E$ the set of bridges and by
 $E_l\subset E$ the set of loops.

\begin{prop}\label{prop_graphic_and_cographic_US}
${\sF _{\Ga }=( F_{\Ga }  , E\setminus E_b)}$ and ${\sG _{\Ga }=(
G_{\Ga } , E\setminus E_l)}$ are unimodular systems. The complexity
of each of the systems is equal to the complexity of the
graph\footnote{Recall that the complexity of a graph is, by
definition, the number of its spanning trees.} $\Ga $.
\end{prop}
Denote by $M_F$ the lattice generated by edges from $E\setminus E_b$
in $F_{\Ga }^*$ and by  $M_G$ the lattice generated by edges from
$E\setminus E_l$ in $F_{\Ga }^*$. Any edge is integral function on $
F_{\Ga }(\Z )$ (respectively on $G_{\Ga }(\Z ) $) and therefore
 $ M_F \subset  F_{\Ga }(\Z )^*$   (respectively
 $ M_G \subset   G_{\Ga }(\Z )^* $). But for
 any spanning tree $T$ in $\Ga $ all the edges $e\notin T$ form a dual basis
 to the basis (\ref{basis_ga}) of $ F_{\Ga }(\Z ) $.
 (Respectively
 all the edges $e\in T$ form a dual basis
 to the basis (\ref{basis_de}) of $ G_{\Ga } (\Z )$.)
Hence $ M_F = F_{\Ga }(\Z ) ^*$ (respectively
 $ M_G =F_{\Ga }(\Z ) ^*$).

To complete the proof, it is necessary to show that any basis of
$M_F$ or $M_G$ consisting of certain edges comes from a spanning
tree in $\Ga $. Consider a set of edges $B = \{ e_1, \ldots , e
_{m-1} \}$ that are not the edges of any spanning tree of $\Ga $.
Denote by $\Ga _1 $ the subgraph of $\Ga $ consisting of all the
edges of $B$ and all the vertices of $\Ga $. Since $\Ga _1 $ is not
a spanning tree it is disconnected and at least one connected
component of $ \Ga _1 $ contains a cycle. Then the set of edges
$E\setminus B$ is not a basis of $F_{\Ga } $ because all the edges
in $E\setminus B$ vanish on that cycle. Similarly  $B$ is not a
basis of $G_{\Ga }  $ because all the edges in $E $ vanish on $\p ^*
(\sum _{v\in V( \Ga _1 ') } v )$ where $ \Ga _1 '$ is one of the
connected components of $ \Ga _1 $. $\square $

The  unimodular systems ${\sF _{\Ga }=( F_{\Ga }  , E\setminus
E_b)}$ and ${\sG _{\Ga }=( G_{\Ga } , E\setminus E_l)}$  are  named
respectively  {\it graphic} and {\it cographic}.

The cographic  unimodular system $\sG  _{\Ga }$ may be constructed
from  the incidence matrix $I _{\Ga }$ in a way similar to that
described in Proposition \ref{prop_sys-matr} It is well-known that
the matrix $I_{\Ga }$ is totally unimodular $N\times m$ matrix (see,
for instance, \cite{Schr}, section 19.3) but
 $\rk I_{\Ga }=m-1$.  So in terms of
Proposition \ref{prop_sys-matr} the cographic unimodular system $\sG
_{\Ga }$ corresponds to totally unimodular matrix obtained from $
I_{\Ga } $ by removal of one arbitrary column.

In the case of planar graphs, the connection between these two
systems is clarified by consideration of the dual graph: the graphic
system of the graph
  $\Ga $ turns out to be the cographic system of the dual graph
   $\widehat{\Ga }$ and vice versa.

\begin{prop}\label{prop_planar_graphs}
Let $\Ga $ be a planar graph embedded in a sphere $S^2$, and
  $\widehat{\Ga }$ be the dual embeded graph.
  Then the graphic
system of the graph
  $\Ga $ is isomorphic to the cographic system of the dual graph
   $\widehat{\Ga }$ and  the cographic system  of the graph
  $\Ga $ is isomorphic to the  graphic
system of the  dual graph
   $\widehat{\Ga }$.
\end{prop}
Embedded graph  $\Ga $ defines a cellular decomposition of the
sphere $S^2$, which provides the chain complex
\begin{equation}\label{complex_for_planar_graph}
    C_2(\Ga ,\R ) \stackrel{\p _2 (\Ga ) }{\to }    C_1(\Ga ,\R )
    \stackrel{\p _1 (\Ga )}{\to }  C_0(\Ga ,\R ).
\end{equation}

(Here
\begin{equation}\label{2-chains}
 C_2 (\Ga ,\R )
=\{ \sum _{\mbox{\tiny faces } f \ \mbox{\tiny of the embeded graph  } \Gamma }
 c_f \cdot f, \quad c_f\in \R \} ,
\end{equation}
orientation of all the 2-cells is induced from the orientation of
$S^2$, and $\p _2 (\Ga )$ is the standard border mapping.)

Since $H_1(S^2)=0$ the sequence (\ref{complex_for_planar_graph}) is
exact, i.e.
\begin{equation}\label{Ker_p1_=_Im_p2}
    \Ker (\p _1 (\Ga ))=\Im (\p _2 (\Ga )).
\end{equation}

We may fix
 the scalar product on $C_2 (\Ga ,\R )$ in the same way as it was
done on $C_0 (\Ga ,\R )$ and $C_1 (\Ga ,\R )$ (the 2-cells form an
orthonormal basis) and identify    $C_2 (\Ga ,\R )$  with its dual.
Then we get the dual sequence to (\ref{complex_for_planar_graph})
\begin{equation}\label{dual_complex_for_planar_graph}
    C_0(\Ga ,\R ) \stackrel{\p _1 (\Ga ) ^*}{\to }    C_1(\Ga ,\R )
    \stackrel{\p _2 (\Ga ) ^*}{\to }  C_2(\Ga ,\R )
\end{equation}
which is nothing else but the sequence
(\ref{complex_for_planar_graph}) for the dual graph $
\widehat{\Ga}$:
\begin{equation}\label{complex_for_dual_planar_graph}
    C_2(\widehat{\Ga} ,\R ) \stackrel{\p _2  (\widehat{\Ga} ) }{\to }
       C_1(\widehat{\Ga} ,\R )
    \stackrel{ \p _1  (\widehat{\Ga} ) }{\to }  C_0(\widehat{\Ga} ,\R )
\end{equation}
and therefore $\p _2  (\Ga )  =\p _1  (\widehat{\Ga} )^* $ and
\begin{equation}\label{Ker_p2=Im_p0*}
    \Ker (\p _1  (\Ga ) )=\Im (  \p _1  (\widehat{\Ga} )^*).
\end{equation}
This implies  that the graphic system for $\Ga $ is the
 cographic system for $\widehat{\Ga }$ and vice versa.$\square $

All  unimodular systems from our previous examples are graphic or
cographic. Unimodular  system defined by matrix (\ref{n=2_N=3}) is
graphic system for theta-graph (two vertices connected by three
edges), or a cographic system for triangle (complete graph with  3
vertices).
 Unimodular  system defined by matrix (\ref{n=2_N=2}) is
graphic system for graph with one vertex and two loops,   or
cographic system for tree with  3 vertices. Unimodular  system from
Example \ref{ex_1N} is cographic system for generalized theta-graph
(two vertices connected by $N$ edges)  or graphic system for polygon
with $N$ vertices.

\section{The lattice.}
Let $\Om =(U,\Xi )$ be an unimodular system, $\Xi =[ \xi _1, \ldots
, \xi _N]$, $\xi _i \in U^*$, and  $ \xi _1, \ldots , \xi _N$
generate over $\Z$ a free abelian group $M\subset U^*$ of rank $n$.
Consider the linear mapping
\begin{equation}\label{def_Phi}
 \Phi _{\Om }: U \to \R ^N \ \mbox{   defined by  } \ \Phi _{\Om }(u)=(\xi _1(u), \ldots ,\xi
_N(u)) .
\end{equation}

Put
\begin{equation}\label{def_V_Om}
 W_ { \Omega } =\Phi _{\Om }(U) \ \mbox{ and } \ L_ {\Omega }=W_ { \Omega }\cap \Z^N
\end{equation}
 where $\Z^N \subset \R^N$
 is the standard lattice consisting of vectors with  integer coordinates.

\begin{prop}\label{Phi_inject}
    The mapping $\Phi _{\Om }$ is injective (i.e. $\dim W_{\Omega }=n$) and
      ${L_ {\Omega }= \Phi _{\Om }(M^*)}$.
\end{prop}
\begin{cor}\label{full_rank} The
     rank of  the lattice
     $L_ {\Omega }$ is equal to   $n$.
\end{cor}

Since  $\Xi $ generates $U^*$ the mapping $\Phi _{\Om }$ is
injective. Choose a basis of M in $\Xi $, say $\xi _1, \ldots , \xi
_n$. (We may change the numeration of elements of
 $\Xi $ if necessary.) Denote by $u_1, \ldots , u_n$ the dual
 basis in $U$. Then the vectors
 \begin{equation}\label{w_i=Fi(u_i)}
  w_1=\Phi _{\Om }(u_1),   \ldots , w_n=\Phi _{\Om }(u_n)
 \end{equation}
 are integer and
 linearly independent in $W_{\Om }$.

 Consider any vector $w\in L_{\Om }$ then $w=\Phi _{\Om }(u)$ for
 certain $u\in U$, then $\xi _i (u)$
 is integer for any $i$, $ 1\le i \le N$. Thus $u\in M^*$. $\square $
\begin{rmkk}\label{remark_matrix_of_Phi}
Denote by $h_1, \ldots , h_N$ the standard basis in $\R ^N$ (i.e.
 the $ k$-th coordinate of $h_k$ is 1, and the rest are zero) then
the matrix of the linear mapping $\Phi _{\Om }$ with respect to the
basis $u_1, \ldots , u_n$ is exactly the totally unimodular matrix
$A_{\Om }$ from (\ref{A_Om}).
\end{rmkk}

Let us fix the standard scalar product on $\R ^N$: $\langle x, y
\rangle = \sum _{i=1}^{N}x_i y_i $ ($x=(x_1, \ldots , x_N)$,
$y=(y_1, \ldots , y_N)$, $x,y \in \R ^N$).

\begin{prop}\label{Kirh}
    The discriminant of  the lattice
     $L_ {\Omega }$  is equal to   the complexity of the unimodular system  $c(\Omega )$.
\end{prop}

The proof repeats a well-known argument proving the classical
Kirhhoff's theorem (\cite{Kirch}): the bases of $U^*$ consisting of
vectors from $\Xi $ are in one-to-one correspondence with bases of
the matrix $ A_{\Om } $ and all these non-zero minors are equal to
$\pm 1$. Therefore the number of such non-zero minors by
Cauchy-Binet theorem is equal  to $\det ( A_{\Om }   ^t  A_{\Om }
)$. And the matrix $  A_{\Om }   ^t  A_{\Om }  $ is nothing but the
Gram matrix of the basis (\ref{w_i=Fi(u_i)}) of the lattice $L_{\Om
}$. $ \square $

We can also define  scalar product on $U$
\begin{equation}\label{scal_prod_on_U}
    \langle u, v \rangle _U= \sum _{i=1}^N \xi (u_i)\xi (v_i),
\end{equation}
then $\Phi _{\Om }$ becomes an isometry of Euclidian spaces $U$ and
$W_{ \Omega }$ and isomorphism of lattices $L_{\Om }$ and $M^*$.
Thus we can identify the spaces $U$ and $W_{\Om }$ by means of the
isomorphism  $\Phi _{\Om }$. The linear functions $\xi _i$ on $U$
may be expressed by
\begin{equation}\label{xi_i(u)=(Fi(u),h_i)}
\xi _i (u)=\langle \Phi _{\Om }(u), h_i \rangle , \qquad u\in U,
\end{equation}
and we shall denote  the corresponding linear functions on $W_{\Om
}$ by the same letter:
\begin{equation}\label{bar_xi_i}
 \xi _i  (w)=\langle w, h_i \rangle , \qquad w\in W_{\Om } ,
\end{equation}
obtaining the same unimodular system $\Om =( W_{\Om }, \Xi )$.

Note that in Proposition \ref{prop_graphic_and_cographic_US} we have
already constructed the lattice $L_{\Om }$ in the case when $\Om $
is  graphic or cographic system for a connected graph $\Ga $. Namely
for the graphic system $\sF _{\Ga }=(F _{\Ga }  , E
 \setminus E_b )$ (we use the notations of Proposition
\ref{prop_graphic_and_cographic_US})
\begin{equation}\label{V_graphic}
W_{\sF  _{\Ga }} =F_{\Ga },
\end{equation}
and
\begin{equation}\label{L_graphic}
L_{\sF  _{\Ga }} =F_{\Ga }(\Z ).
\end{equation}
while for cographic system
 $\sG _{\Ga }=(G _{\Ga }  , E
 \setminus E_l )$
\begin{equation}\label{V_cographic}
W_{\sG _{\Ga }} =G _{\Ga }
\end{equation}
and
\begin{equation}\label{L_cographic}
L_{\sG _{\Ga }} =G _{\Ga } (\Z ).
\end{equation}

In addition to the bases  constructed in the proof of the
 Proposition \ref{prop_graphic_and_cographic_US} the lattice
 $L_{\sG _{\Ga }}$ has a remarkable basis
which is used in the proof of the Kirhhoff's theorem. Fix a vertex
$v_0\in V$, then the  elements (\ref{codiff})
\begin{equation}\label{cut_defined by_vertex}
 {\p ^*(v)= \sum _{e\in E\setminus E_l}
  i(e,v)\cdot e}
\end{equation}
 for all $v\in V(\Ga )$,  $v\ne v_0$, form a basis
  of   $L_{\sG _{\Ga }}$.
The Gram matrix of this basis is obtained from the Laplacian matrix
$I_{\Ga }^T I_{\Ga }$ by deleting the row and the column
corresponding to the vertex $v_0$. Using Proposition
 \ref{prop_graphic_and_cographic_US} we get the Kirhhoff's theorem:
the determinant of this Gram matrix   is equal to   the complexity
of the graph $\Ga $.
 \begin{cor}\label{cor_A_N_*_for_complete_graph}
 The lattice  $L_{\sG  _{ K_N }}$ for the
 complete graph    $K_N$ on $N$ vertices is $
  A_{N-1}^*$ (see \cite{Conv})
--- integer lattice congruent to the lattice
dual to the root lattice of type  $A_{N-1}$.
 \end{cor}
Indeed, the Gram matrix for $K_N$ is
\begin{equation}\label{Gram_matrix_for_complete_graph}
   \left(
      \begin{array}{ccccc}
        N-1 & -1 & -1 & \ldots  & -1 \\
        -1 & N-1 & -1 &  \ldots  & -1 \\
        -1 & -1 & N-1 &  \ldots  & -1 \\
         \ldots  &  \ldots  &  \ldots  &  \ldots  &  \ldots  \\
        -1 & -1 & -1 &  \ldots  & N-1 \\
      \end{array}
    \right) .
\end{equation}
which is exactly the the Gram matrix for  $ A_{N-1}^*$ (see
\cite{Conv}, chapter 4, (77)). $\square $

Note that the complexity of the complete graph is the number of
 of trees on
$N$ labeled vertices and the determinant of
(\ref{Gram_matrix_for_complete_graph}) is $N^{N-2}$, so Propositions
 \ref{prop_graphic_and_cographic_US} and \ref{Kirh} provide the well-known proof
 of Cayley's formula.

Now let's continue the examples for for generalized theta-graph.

\begin{ex}\label{ex_gen_theta_graph1}
We have already mentioned that  cographic system for generalized
theta-graph $\Th  _N $ (two vertices connected by $N$ edges $e_1,
\ldots , e_N$) or graphic system for dual graph  $\widehat{\Th  _N
}$ (which is a polygon with $N$ vertices) is the system $\Si _N$
from Example \ref{ex_1N}. In this case $W_{\Si _N }$ is
one-dimensional subspace of ${C_1 (\Th  _N , \R )=C_1 (\widehat{\Th
_N }, \R ) =\R ^N}$, ${W_{\Si _N } =\{ (x_1,\ldots , x_N)\in \R ^N \
, \ x_1=x_2= \ldots =x_N \}}$, and $L_{\Si _N  }$ consists of such
vectors with integer entries. The discriminant of the lattice $L_{
\Si _N  }$ is $N$, which is also the complexity of the unimodular
system $c( \Si _N  )$ and the complexity of each of  the graphs $\Th
_N $ and  $\widehat{\Th  _N }$.
\end{ex}
\begin{ex}\label{ex_gen_theta_graph2}
Consider the graphic system  $\La _N =(\Ker \p , \{ e_1, \ldots ,
e_N \})$ for the generalized theta-graph $\Th  _N $ from the
previous example. (Of course $\La _N$ is the cographic system for
$\widehat{\Th _N}$ as well.) In this case $\Ker \p $ is $N-1$
dimensional subspace of $C_1 (\Th _N , \R )=\R ^N$, so
 $ { W _{\La _N} =\{ (x_1,\ldots , x_N)\in \R ^N \ , \ \sum
_{i=1}^N x_i =0 \}}$ and $  L _{\La _N}$ consisting of such vectors
with integer entries is the root lattice of type $A_{N-1}$. The
discriminant of this lattice is $N$, which is also the complexity of
the unimodular system $c( \La _N  )$ and the complexity of each of
the graphs $\Th _N $ and  $\widehat{\Th  _N }$.
 Roots of
this lattice are the cycles of length 2 in $\Th  _N $ (or the 2-cuts
of the polygon $\widehat{\Th  _N }$).
\end{ex}

\section{Direct sums of unimodular systems.}
Consider two unimodular systems $\Omega  =( U, \Xi )$ and $\Omega '
=( U', \Xi ')$. Standard inclusions $U^* \subset (U\oplus U')^*$ and
$U'^* \subset (U\oplus U')^*$ alow to consider $ \Xi \cup \Xi '$ as
a unimodular system in  $ (U\oplus U')^*$. Then we have the direct
sum of mappings $\Phi _{\Om \oplus \Om '} = \Phi _{\Om } \oplus \Phi
_{\Om '}: U \oplus U' \to \R^N \oplus \R^{N'}$ and the space $W_{\Om
\oplus \Om '}$ with the scalar product $\langle \cdot , \cdot
\rangle  $ is the orthogonal sum of $W_{\Om }$ and  $W_{\Om '}$ with
the corresponding scalar products.

Note that if the matrices of the mappings $\Phi _{\Om }$ and  $\Phi _{\Om '}$ are
consequently  $ \left(
                       \begin{array}{c}
                         E \\
                         A \\
                       \end{array}
                     \right) $
and $ \left(
                       \begin{array}{c}
                         E \\
                         A' \\
                       \end{array}
                     \right) $
                     (see (\ref{A_Om}))
then the matrix of the mapping $\Phi _{\Om \oplus \Om '}$ is
\begin{equation}\label{matrix_of_direct_sum}
 \left(
                       \begin{array}{cc}
                         E &0 \\
                         0 & E\\
                         A & 0\\
                         0& A' \\
                       \end{array}
                     \right)
\end{equation}

It is clear that compexity is multiplicative:
\begin{equation}\label{mult_compl}
    c(\Om \oplus \Om ')=  c(\Om ) c( \Om ') .
\end{equation}

The most trivial unimodular system consists of one linear form $\xi
\in U^*$ for one-dimensional vector space $U$. In the Example
\ref{ex_Upsilon} we denoted this system by $\Upsilon $. Note that
$c(\Upsilon )=1$ and
\begin{equation}\label{c_Om_plus_Ups}
    c(\Om \oplus \Upsilon )=  c(\Om ) .
\end{equation}
Note that this in the only case when the mapping $\Phi  _{\Omega }$
is surjective (and therefore bijective):
\begin{prop}\label{bijective_Phi}
$\Phi _{\Omega }$ is bijective if and only if $\Omega =  \Upsilon
^m$ for some $m$. $ \square $
\end{prop}

It is easy to characterize  unimodular  systems containing $\Upsilon
$ as a direct summand. Let  $\Omega =  \Upsilon  \oplus \Omega '$,
where $ \Omega ' = (U' , \Xi ' ) $  and   $\Upsilon = (\tilde{U}, \{
\tilde{\xi} \} )$, $\dim \tilde{U} =1$, $\tilde{u}$  is the
generator of
 $\tilde{U}$ and $\tilde{\xi}$ the dual generator of $ \tilde{U}^*
 $. Then the matrix (\ref{matrix_of_direct_sum}) takes the following
 form:
\begin{equation}\label{matrix_of_Om_plus_Ypsilon}
 \left(
                       \begin{array}{cc}
                         1 &0 \\
                         0 & E\\
                         0 & 0\\
                         0& A' \\
                       \end{array}
                     \right)
\end{equation}
and $\Phi _{\Om } ( \tilde{u}, 0)=(1,0,0,\ldots ,0)$ --- the first
vector of the standard basis $h_1, \ldots , h_N$ of $\R ^N$(see
Remark \ref{remark_matrix_of_Phi}).  So $h_1\in \Phi _{\Om }
(U)=W_{\Om }$ and $|h_1|=1$. Moreover, any vector of unit length in
$L_{\Om }$ is $\pm h_i$ and therefore comes (after suitable
renumbering) from the above construction.

\begin{prop}\label{one_dim_summand}
    $\Omega = \Omega ' \oplus  \Upsilon $ if and only if $L_{\Omega }$ contains
    a vector $w$ such that $|w|=1$. $
\square $
\end{prop}

The same in terms of matrix $A_{\Om }$ (\ref{A_Om}):
\begin{prop}\label{prop_summand_Ups_matrix}
 $\Omega = \Omega ' \oplus  \Upsilon $ if and only if there is a
 column of the matrix $A_{\Om }$ having only one nonzero entry.$
\square $
\end{prop}

\begin{prop}\label{prop_eliminate_Upsilon}
For any unimodular system $\Om $ there is a canonical representation
\begin{equation}\label{eliminate_Upsilon}
\Om =\Om '\oplus \Upsilon ^s
\end{equation}
  where the system $\Om '$ has
no direct summand $\Upsilon $.
\end{prop}
The decomposition $U=U' \oplus U_{\Upsilon ^s}$ may be defined as
follows: consider the subspace in $U$ generated by all the vectors
in $u\in M$ such that $\langle u, u \rangle _U=1$ and denote its
dimension by $s$; then it has exactly $s$ such vectors (up to sign)
corresponding to $s$ elements of $\Xi $. Each of these elements
defines $\Upsilon $ so  this $s$-dimensional space is exactly  $
U_{\Upsilon ^s} $. Let $U'$ be the orthogonal complement to   $
U_{\Upsilon ^s} $  in $U$; then the remaining
 $N'=N-s$ elements of $ \ \Xi \ $
 form the unimodular system $\Om '$ on $U'$. Note that there is
an orthogonal decomposition  $\R ^N _{\Om } = \R ^{N'}_{\Om '}
\oplus \Phi_{\Om } (U_{\Upsilon ^s})$ and $ \left. \Phi _{\Om '} =
\Phi _{\Om } \right| _{U'}$. $ \square $

\begin{prop}\label{prop_graph_eliminate_Upsilon}
For graphic and cographic systems for graph $\Ga $ direct summands $
\Upsilon $ correspond respectively to loops and bridges, and
representation (\ref{eliminate_Upsilon}) can be described as
follows.
\begin{equation}\label{Ups_graphic}
     \sF  _{\Ga }
    = \sF  _{\Ga ' }  \oplus \Upsilon ^{s' }
\end{equation}
where $\Ga ' $ is obtained from  $\Ga  $ by deleting all the loops,
and $s' $ is the number of loops in  $\Ga  $.
\begin{equation}\label{Ups_cographic}
     \sG  _{\Ga }
    =  \sG  _{\Ga ''}  \oplus \Upsilon ^{s''}
\end{equation}
where $\Ga ''$ is obtained from  $\Ga  $ by contracting all the
bridges, and $s''$ is the number of bridges in  $\Ga  $.
\end{prop}

\section{The dual system.}\label{section_duality}
In this section we construct the Gale dual unimodular system $\Omega
^{\bot }$. Let $\Om =(U, \Xi)$, consider the data (\ref{def_Phi})
and (\ref{def_V_Om}). Take the orthogonal complement to $W_{\Omega
}$ in $\R ^N$, namely
\begin{equation}\label{V^bot}
W^{\bot }_{\Omega } = \{ z\in \R ^N , \ \langle w, z\rangle =0 \
\forall w \in W_{\Omega } \} .
\end{equation}
 Each vector $h_i$ from the standard basis of $\R ^N$
 (see Remark \ref{remark_matrix_of_Phi}) provides a linear function
$\eta  _i (z)= \langle h_i , z \rangle $ on $W^{\bot }_{\Omega }
\subset \R ^N $. Of course the function $\eta _i $  may not be
non-zero on $ W_{\Omega }^{\bot}$. Suppose that some  $\eta _i$, say
$\eta _1$, vanishes on $ W_{\Omega }^{\bot}$. This means that $
\langle h_1 , z \rangle =0$ for all $z \in W^{\bot }_{\Omega }$ and
therefore $h_1\in W_{\Omega }$, so  $h_1 \in L _{\Omega } $ and this
is just the case described in Propositions \ref{one_dim_summand} and
\ref{prop_summand_Ups_matrix}.

\begin{prop}\label{zero_funct_on_dual}
$\eta  _i $ vanishes on $V^{\bot }$ if and only if $\xi _i$
generates a direct summand of $\Om $ isomorphic to $\Upsilon $. $
\square $
\end{prop}
\begin{prop}\label{duality}
Let    $\Omega = \Omega ' \oplus  \Upsilon ^s$, where $\Omega '=(U',
\Xi ')$ has no direct summand $\Upsilon $. Then $V _{\Om } ^{\bot }=
V _{\Om '} ^{\bot }$. $ \square $
\end{prop}
We have seen that $\eta _i $ vanishes on $W_{\Om } ^{\bot }$ if and
only if $\xi _i$ corresponds to a trivial direct summand $\Upsilon
$. Let us collect al  non-zero functions  $\eta _i $ and put
\begin{equation}\label{Xi^bot}
\Xi ^{\bot } = [ \eta _i, \quad  \xi _i \in \Xi ' ].
\end{equation}

\begin{prop}\label{def_of_dual_system}
Put  $\Omega ^{\bot} =( W _{\Om }^{\bot },  \Xi  ^{\bot } )$. Then
         $\Omega ^{\bot}$  is a unimodular system of
         linear functions on $W_{\Om }^{\bot }$.
\end{prop}
It is sufficient to prove this statement for the case when $\Omega $
has no direct summand $\Upsilon $. Let $|\Xi |=N$, $\dim  W _{\Om }
=n$, $\dim  W _{\Om } ^{\bot } =m$, where $m+n=N$. According to
Definition \ref{def_unumod_sys} we need to prove that all linearly
independent subsets of $\Xi ^{\bot } $ define over $\Z $ the same
abelian subgroup of $(W _{\Om }^{\bot })^*$. The steps of the proof
are collected in the following lemma.
\begin{lem}\label{lem_dop_bazisov} Let   $\Omega   =( U,  \Xi   )$ be
an unimodular system having no direct summands $\Upsilon $ and $(W
_{\Om }^{\bot })^*$ and  $\Xi ^{\bot } $ defined by (\ref{V^bot})
and (\ref{Xi^bot}). Then \\
a) Linear functions
  $\eta _{j_1}, \ldots \eta _{j_{N-n}}$ are linearly independent if
  and only if $\xi _{i_1}, \ldots \xi _{ i_n }$ are linearly
  independent,
where $\{i_1 , \ldots   i_n \} \sqcup \{j_1 , \ldots   j_{N-n} \} =
\{
1, \ldots , N \}$.\\
b) The lattice  $W_{\Om }^{\bot} \cap \Z ^N$ has rank $N-n$.\\
c) Linear functions $\eta _1, \ldots ,\eta _N$ generate over $\Z $
the group $( W_{\Om }^{\bot} \cap \Z ^N)^*\subset ( W_{\Om } ^{\bot})^*$.\\
d) Any  linearly
  independent $\eta _{j_1}, \ldots \eta _{j_{N-n}}$
  form a basis of  $(W_{\Om } ^{\bot}
\cap \Z ^N)^* $  over $\Z $.
\end{lem}

In a) we characterize  linearly independent subsets of $\Xi ^{\bot }
$, in b) and c) we describe the abelian group generated by $\eta _1,
\ldots ,\eta _N$   over $\Z $, and then d) is just the statement we
need to prove. $ \square $

{\it Proof of Lemma \ref{lem_dop_bazisov} }. For simplicity we may
renumber the elements of $\Xi $ so that  $\{i_1 , \ldots   i_n \} =
\{ 1, \ldots , n \}$ and  $ \{j_1 , \ldots   j_{N-n} \} =\{ n+1,
\ldots , N \}$.

 Let's assume that
 $\xi _1, \ldots , \xi _n$ are not linearly independent. Then
 there is a nontrivial relation
$a_1 \xi _1+ \ldots +a_n \xi _n =0$. Therefore ${a = (a_1, \ldots ,
a_n, 0, \ldots , 0) \in W_{\Om }^{\bot }}$ is nonzero, but $\eta _j
(a)=0$ for $j=n+1, \ldots , N$ and hence $\eta _{n+1} , \ldots ,
\eta _N$
 is not a basis of $(W_{\Om }^{\bot })^*$.

Next let's assume that
 $\xi _1, \ldots , \xi _n$ are linearly independent, then
they form a basis of $U^*$. Let $u_1, \ldots , u_n$ be the dual
basis of $U$; put $z_i=\Phi _{\Om }(u_i)$. Then
\begin{equation}\label{v_i=}
z_i  = h_i + \sum _{q>n} a_{q,i} h_q, \quad 1\le i \le n,
\end{equation}
 where
$a_{q,i}$ are the entries of the matrix $A_{\Om }$ from (\ref{A_Om}
). (Recall that $A_{\Om }$ is the matrix of the mapping $\Phi _{\Om
}$, see remark \ref{remark_matrix_of_Phi}.)  Therefore for $j>n \ $
$ \xi _j = \sum _{p\le n} a_{j,p} \xi _p$. Put
\begin{equation}\label{w_j}
w_j = -h _j + \sum _{p\le n} a_{j,p} h _p \  \mbox{  for  } j>n,
\end{equation}
 then
$\langle z_i , w_j \rangle =0$ for $i\le n$, $j>n$, so $w_j \in
W_{\Om }^{\bot }$. The elements $w_j$ for $j>n$ are linearly
independent and therefore form a basis of $W_{\Om }^{\bot} $ over
$\R $, which proves item a). Note that these elements are integer:
$w_j \in W_{\Om }^{\bot} \cap \Z ^N$, therefore rank of the lattice
$ W_{\Om }^{\bot} \cap \Z ^N$ is $N-n$, which proves item b).

Let's prove that $w_{n+1}, \ldots , w_N$ also form a basis of the
lattice $W_{\Om }^{\bot} \cap \Z ^N$. Let $w=\sum _{i=1}^N x_i h_i$
be an integer element of $W_{\Om }^{\bot }$, i.e. all the $x_i$ are
integer. Consider $\tilde w = w + \sum _{j>n} x_j w_j$, then $\tilde
w \in W_{\Om }^{\bot }$. But for $j>n \ $ $\langle \tilde w ,
e_j\rangle =0 $ and therefore $\tilde w = \sum _{i\le n} \tilde x_i
h_i$. Then for $i\le n \ $ $ \langle \tilde w , z_i\rangle =\tilde x
_i $ and since $ \tilde w \in W_{\Om }^{\bot}$ all $\tilde x _i =0$.
Therefore $\tilde w = 0$ and ${w=-\sum _{j>n} x_j w_j}$ is the
expansion of $w$ over $w_{n+1}, \ldots , w_N$ with integer
coefficients, so  $w_{n+1}, \ldots , w_N$ is the basis of the
lattice $W_{\Om }^{\bot} \cap \Z ^N$, which proves c). Now it only
remains to note that $\eta _k (w_j)=\langle h _k , w_j\rangle =\de
_{k,j} $ and therefore
  $\eta _{n+1}, \ldots , \eta _N$ is the  basis of  the
  lattice  $(W_{\Om }^{\bot} \cap \Z ^N)^*$, dual to the  basis
  $w_{n+1}, \ldots , w_N$. So we see that any maximal
   linearly independent
subset in $\Xi ^{\bot }$ generates over $\Z $ the same  lattice
${(W_{\Om }^{\bot} \cap \Z ^N)^*}$, which proves d). $\square $

Note that if we choose the basis $\eta _{n+1}, \ldots , \eta _N$ of
the unimodular system $\Om ^{\bot }$ then the matrix $A_{\Om ^{\bot
}} $ defined in  (\ref{A_Om}) is just the matrix
\begin{equation}\label{A_Om^bot}
   A_{\Om ^{\bot }} =
   \left(
           \begin{array}{c}
          \widetilde{A}^T \\
             -E  \\
           \end{array}
         \right) ,
\end{equation}
where $E$ is the unit $(N-n)\times (N-n)$-matrix and $
\widetilde{A}^T$ is the transposed matrix   $\widetilde{A}$ from
(\ref{A_Om}). Recall that (\ref{w_j}) corresponds to the case when
$\Om $ has no direct summands $\Upsilon $. In general case a direct
summand $\Upsilon $ correspond to a column of $A_{\Om }$ with only
one nonzero element, which implies zero column of $ \widetilde{A}$.
This leads to the following matrix description of the matrix
corresponding to $\Om ^{\bot }$.
\begin{cor}\label{cor_matrix_of_Om_bot}
The dual unimodular system $\Om ^{\bot }$ is isomorphic to the
unimodular  system of nonzero rows of the matrix (\ref{A_Om^bot}).
\end{cor}
\begin{cor}\label{cor_Seymour_autodual}  The Bixby-Seymour unimodular
system from Example \ref{ex_Seymour} is isomorphic to its dual.
\end{cor}
Indeed, the matrix transposed to the matrix $\widetilde{A}$ (five
bottom rows of (\ref{Seymour})) can be obtained from $\widetilde{A}$
by a permutation of its rows. $\square $
\begin{cor}\label{cor_c(Om)=c(Om_bot)}
\begin{equation}\label{c(Om)=c(Om_bot)}
     c(\Omega )=c(\Omega ^{\bot}  ) .
\end{equation}
\end{cor}
Lemma  \ref{lem_dop_bazisov}a)  provides one-to-one correspondence
between   bases of $\Xi $ and  bases of  $\Xi ^{\bot}$ which implies
that the  complexities are equal. $\square $
\begin{prop}\label{prop_no_Upsilon_in_dual}
   The Gale dual system $\Omega ^{\bot}$ has no direct summands
   $\Upsilon $.
\end{prop}
Assume that $\Omega $ has a direct summands
   $\Upsilon $. Then there is certain $\tilde{w}  \in V^{\bot}$ such that
  $\eta _j (\tilde{w}) \ne 0$ for  only one $j$, say for  $j=1$ and (up
to scalar factor)   $\eta _1 (\tilde{w}) =1$. Since $ \eta _j
(\tilde{w}) = \langle h _j,\tilde{w} \rangle $ and $\tilde{w}=\sum
\langle h _i,\tilde{w} \rangle h_i$ it means that $\tilde{w}=h_1$,
so $h_1 \in V^{\bot}$. Therefore $\xi _1 $ is zero function on $U$
(see (\ref{xi_i(u)=(Fi(u),h_i)})), which contradicts the Definition
of unimodular system. $\square $
\begin{prop}\label{pror_double_dual} Let   $\Omega   =( U,  \Xi   )$ be
an unimodular system having no direct summands $\Upsilon $. Then
$(\Omega ^{\bot}) ^{\bot} =\Om $.
\end{prop}
Since $\R ^N =W_{\Om } \oplus W_{\Om }^{\bot} \ $ $(W_{\Om
}^{\bot})^{\bot} =W_{\Om }$ and all
 $\eta _j $ are nonzero functions on $ W_{\Om }^{\bot}$,
 so $(\Xi ^{\bot})^{\bot} =\Xi $. $\square $

 Now it is clear that in general case
 double dual system returns the initial one
omitting  direct summands $\Upsilon  $.
\begin{prop}\label{prop_duality1}
Let    $\Omega = \Omega ' \oplus  \Upsilon ^s$ where $\Omega '$ has
no direct summand $\Upsilon $. Then $(\Omega ^{\bot})^{\bot } =
\Omega '$.  $\square $
 \end{prop}

Note that a direct summand $\Upsilon $ in a graphic unimodular
system corresponds to a loop in the graph, and  a direct summand
$\Upsilon $ in a cographic unimodular system corresponds to a
 bridge in the graph. So for a graph $\Gamma $ without
 loops and bridges all the edges are nonzero linear functions on
 $ F _{\Ga }$  and $  G _{\Ga }$ from the orthogonal
decomposition (\ref{C_1=OplusFG}) and therefore unimodular systems
 $\sF _{\Ga } $ and $\sG _{\Ga }   $ are Gale dual to each other.
 \begin{prop}\label{Gale_duality_for_graphs} A graph $\Gamma $ has
 no loops and bridges if
 and only if it's graphic and cographic systems
 $ F _{\Ga }$  and $  G _{\Ga }$  are Gale dual to each other. $\square $
   \end{prop}
For a  graph having certain  loops and bridgs we use Propositions
\ref{prop_duality1} and
  \ref{prop_graph_eliminate_Upsilon} and conclude the following.
\begin{prop}\label{stabilization_of_graph} Suppose that graph $\overline{\Gamma }$ is obtained
from a graph $\Gamma $ by contracting all its bridges and removal of
all its loops. Then the graphic (cographic) system of
$\overline{\Gamma }$ is the Gale dual of the graphic (cographic)
system of $\Gamma $. $\square $
\end{prop}

We have seen in Examples \ref{ex_gen_theta_graph1} and
\ref{ex_gen_theta_graph2} that the graphic system for the theta
graph $\Theta _N$ is a root lattice (namely $A_{N-1}$) and the
cographic system has multiple elements (namely one element repeated
$ N $ times). It turns out that the same is true in the general
case: the presence of roots in the lattice $L_{\Om }$ for some
unimodular system $\Om $ is always associated with the presence of
multiple elements in the dual system  $\Om ^{\bot }$ and vice versa.

 Consider an
unimodular system  $\Om =(U, \Xi )$, $\Xi  =[ \xi _1, \ldots , \xi
_N]$, having multiple elements, say $ \xi _1= \xi _2$. Then the
vector $w=(1,-1,0,0,0, \ldots ,0)\in \R^N$ is orthogonal to $\Phi
_{\Om }(U)$ and therefore $w\in W_{\Om }^{\bot}$. So $w$ is a root
in  the dual lattice $L_{\Om }^{\bot}$.

Note that if $\xi _1 = \xi _2 = \ldots = \xi _m $ then the
corresponding roots in  $L_{\Om }^{\bot}$ generate a sublattice of
$A_{m-1}$ type.

Conversely, it is not hard to see that any root $w \in L_{\Om
^{\bot} }$ corresponds to a pair of equal (up to sign change)
elements in $\Xi $. Since $w\in \Z ^N$ and $|w|^2=2$ the vector $w$
has exactly two nonzero coordinates, say the first and the second.
We may choose the signs of $\xi _1 $ and $\xi _2$ so that $w=
(1,-1,0,0,0, \ldots ,0)$. Since $w$ is orthogonal to $W_{\Om }$ then
for any $z\in W_{\Om }$ $\langle z, w \rangle =\xi _1 (z) - \xi _2
(z)=0$ and therefore $\xi _1 = \xi _2$ on $W_{\Om }$. So we get the
following statement.

 \begin{prop}\label{roots_and_mult} A unimodular system $\Om$ has multiple
 elements $\xi _1 = \ldots =\xi _m$
 if and only if
 the lattice of the dual system $ L(\Om ^{\bot })$ has a root
 sublattice of $A_{m-1}$ type. $\square $
\end{prop}
\begin{cor}\label{cor_roots_and_mult} An unimodular system $\Om $ has no multiple elements if and only if the
lattice of the dual system $L(\Om ^{\bot })$ has no roots. $\square
$
\end{cor}

\section{The polytope.}\label{section_polytope}
Consider the polytope $\Delta _{\Om } \subset W_{\Om }$ defined by
inequalities $|  \xi _i (w)| \le 1$, $i=1, \ldots , N$. This means
that $ \Delta _{\Om } $ is the intersection of $W_{\Omega }$ and the
standard cube defined by inequalities  $|x_i |\le 1$, $i=1, \ldots ,
N$ in $\R ^N$.
\begin{prop}\label{polytope}
$ \Delta  _{\Om }$ is a reflexive lattice polytope in $W_{\Om }$.
\end{prop}
Any vertex of $\De _{\Omega }$ is defined by $n$ equations $  \xi
_{i_1}= \pm 1 ,$ \ldots  $,  \xi _{i_n}= \pm 1$ for some linearly
independent functions $\xi _{i_1} , \ldots , \xi _{i_n}$ (and
certain choice of signs). But since these functions form a basis of
$L _{\Om } ^*=M$ the solution is integer, so all the vertices of
$\De _{\Om }$ are integer. All the facets of $\De _{\Om } $ are
given by integer equations  $\xi _{i} (u)= \pm 1$ so $ \De _{\Om }$
is reflexive.$ \square $

By definition the polytope $\De _{\Omega }$ is centrally symmetric.
Let us prove that any $\xi \in \Xi $ defines a pair of parallel
facets of $\De _{\Om }$.  Suppose the opposite for say  $\xi _1$,
then the inequality
\begin{equation}\label{xi_le_1}
     \xi _1 (w) \le 1
\end{equation}
  is the consequence of certain other  inequalities
\begin{equation}\label{xi_i_le_1}
     \xi _{k_i} (w) \le 1
\end{equation}
  for some $\xi _{k_1}, \ldots , \xi _{k_m}$ from $\Xi $. Then by the
affine form of the Farkas lemma (see \cite{Schr}, Corollary 7.1h)
the inequality (\ref{xi_le_1}) is a nonnegative linear combination
of inequalities (\ref{xi_i_le_1}) i.e there exist some
 nonnegative $\al _1, \ldots , \al _m$ such that
\begin{equation}\label{lin_comb}
\xi_1
=\sum \al _i \xi _{k_i}
\end{equation}
(left side of the inequality) and
\begin{equation}\label{sum_alpha=1}
 1 =\sum \al _i
\end{equation}
(right side of the inequality). Omitting all zero values of $ \al _i
$ we may assume that all   $ \al _i $ are positive. Fixing
coordinates as above we may consider (\ref{lin_comb}) as a relation
on the rows of the totally unimodular matrix $ A_{\Om }=(a_{p,q})$
(see  (\ref{A_Om})). For the first column we have
\begin{equation}\label{1_column}
 1=\sum \al _i a _{1,k_i}
 \end{equation}
 (since $
a_{1,1}=1$). Comparing with (\ref{sum_alpha=1}) we conclude that all
\begin{equation}\label{Phi_1_ki_=1}
   a _{1,k_i} =1.
 \end{equation}
Suppose that some other $a_{l,k_i} \ne 0$ for $l\ne 1$, then for the
$l$-th column we get a nontrivial relation
\begin{equation}\label{l_column}
 0=\sum \al _i a _{1,k_i}
 \end{equation}
 (since $ a_{l,1}=0$ for $l\ne 1$) and (\ref{l_column})
 has at least one   nonzero term. Then   there are nonzero  terms
 with different signs
 in (\ref{l_column}), say $a_{1,k_s}=1$
 and $a_{1,k_t}=-1$ for some $s$ and $t$.
Taking first and $l$-th column and $s$-th  and $t$-th row we obtain
$2\times 2$ minor
\begin{equation}\label{2minor}
     \left(
       \begin{array}{cc}
         1 & 1 \\
         1 & -1 \\
       \end{array}
     \right)
\end{equation}
which contradicts the total unimodularity of  $A_{\Om }$. So
 all $a_{l,k_i} =0$ for $l>1$ and therefore
 $\xi _{k_1}= \ldots = \xi _{k_m} = \xi _1$.
So this is the only case when inequality  (\ref{xi_le_1}) is a
consequence of other   inequalities defining $\De _{\Om }$ and in
this case all  the inequalities (\ref{xi_le_1}) and
(\ref{xi_i_le_1}) coincide and define the same facet of $\De _{\Om
}$. We have proved the following statement.
\begin{prop}\label{faces_polytope}
Pairs of parallel facets of the polytope  $ \Delta _{\Om }$ are in
one to one correspondence with the set of elements of $\Xi $.
$\square $
\end{prop}

\begin{prop}\label{prop_zonotop}
    Polytope  $ \Delta _{\Om }$ is a zonotope.
\end{prop}
By definition (see \cite{Z}) a zonotope is a projection of a cube.
Our polytope  $ \Delta _{\Om }$ is a section of the standard cube
$[-1;1]^N \subset \R ^N$. Let's prove that  $ \Delta _{\Om }$ is the
image of this cube under the
 projection $\pi : \R ^N \to W_{\Om }$ parallel to the space $W_{\Om
 } ^{\bot }$. It is sufficient to prove that any line $l$ passing
 through a vertex $a=(a_1, \ldots , a_N)$ of  $ \Delta _{\Om }$
 parallel to $W_{\Om
 } ^{\bot }$ is a supporting line for the cube.
Let's assume that $l$ is defined by the vector $y=(y_1, \ldots
y_N)\in
 W_{\Om }^{\bot}$, then all the points of $l$ are given by $a+ty$
 for $t\in \R$. Let's prove that these points can not be internal
 points of the cube. All the coordinates of $a$ are $0$ and $\pm 1$;
 after suitable renumbering of elements of $\Xi $ and suitable  choice of
 signs of $\xi _i$ we may assume that $a_1=a_2=\ldots =a_m=1$ and
 $a_{m+1}=\ldots =a_N=0$. Since $a\in W_{\Om }$ and $y\in W_{\Om
 }^{\bot }$ \ $\sum _{i=1} ^m y_i =0$ and therefore $\exists \ $ $i_1$
 such that $v_{i_1}>0$ and   $\exists \ $ $i_2$
 such that $v_{i_2}<0$. Then  $a_{i_1}+ty_{i_1}>1$ for $t>0$ and
  $a_{i_2}+ty_{i_2}>1$ for $t<0$, so the points  $a+ty$ lie outside
  the cube for all $t\ne 0$. $\square $
\begin{ex}\label{ex_polytop_Ups}
 For the trivial one-dimensional unimodular system $\Upsilon $ (see Example
 \ref{ex_Upsilon}) the polytope $\De _{\Upsilon }$ is the segment
 $[-1;1]$.
\end{ex}
\begin{ex}\label{ex_poly_gen_theta_graph1}
For the  cographic system  $\Si _N$ for generalized theta-graph $\Th
_N $ (see Example \ref{ex_gen_theta_graph1}) $W_{\Si _N }$ is
one-dimensional subspace of \\
 ${C_1 (\Th _N , \R ) =\R ^N}$,
${W_{\Si _N } =\{ (x_1,\ldots , x_N)\in \R ^N \ , \ x_1=x_2= \ldots
=x_N \}}$, and so $\De _{\Si _N }$ is the segment $[(-1,\ldots
,-1);( 1,\ldots , 1)]$.
\end{ex}
\begin{ex}\label{ex_poly_gen_theta_graph2}
For the graphic system  $\La _N $ for $\Th  _N $  (see Example
\ref{ex_gen_theta_graph2}) the lattice $L_{\La _N} $ is the root
lattice of type $A_{N-1}$,
\begin{equation}\label{A_N_lattice}
   L _{\La _N} =\{ (x_1,\ldots , x_N)\in \Z ^N \ , \ \sum
_{i=1}^N x_i =0 \}
\end{equation}
and the polytope $\De_{\La _N}$ is defined by $2N$ inequalities
$|x_i | \le 1$.  Geometrically, it can be  described   as the
intersection of two simplices that are centrally symmetrical to each
other, one of which is given by its vertices $P_1, \ldots , P_N$
where $P_k=(-1,-1, \ldots, -1,  N-1, -1, \ldots, -1)$ (all the
components are $-1$, except for the $k$-th component, which is
$N-1$).
\end{ex}
\begin{ex} \label{poly_K_N}
The
lattice $L_{\Om _{ K_N }^{\mbox{\tiny cographic}}}$ for the
 complete graph    $K_N$ on $N$ vertices is $ N\cdot  A_{N-1}^*$, where  $A_{N-1}^*$
 is the lattice dual to the root lattice of type  $A_{N-1}$
 (see Corollary \ref{cor_A_N_*_for_complete_graph}).

The structure of the polytope $\De _{\sG_{ K_N }}$ may be clarified
using the mapping $\p ^* : C_0(  K_N , \R ) \to  G_{ K_N }$. (Recall
that  $G_{ K_N } =\Im \p ^*$, see (\ref{C_1=OplusFG}).) Let us
denote the vertices of $K_N$ by $v_1, \ldots , v_N$ and the edges by
$e_{i,j}$ (the edge incident to $v_i$ and $v_j$). Recall
(\ref{0-chains}):
 $ C_0 (K_N ,\R ) =\{ \sum _1^N
a_i  v_i, \quad a_i\in \R \} = \R ^N$ and consider the the
hyperplane $H$ defined by the equation $\sum _1 ^N a_i =0$. $H$ is
orthogonal to the kernel of $\p ^*$ (for $\Ker \p ^* $ is
one-dimensional space $ \{ \la  \sum _1^N   v_i \quad , \la \in \R
\} $) and so  $\p ^*$ is the composition of the orthogonal
projection $\pi : C_0 (K_N ,\R )\to H$ and the restriction of $\p
^*$ to $H$. $\pi (v_k) = v_k - \frac 1N \sum _1^N v_i $ so  we can
compare the Gram matrix of the basis $ \pi (v_1), \ldots ,  \pi
(v_{N-1})$ of $H$  and the Gram matrix
(\ref{Gram_matrix_for_complete_graph})  of the basis $ \p ^* (v_1),
\ldots ,  \p ^* (v_{N-1})$ of $\Im \p ^*$ and conclude that
 the restriction of $\p
^*$ to $H$ is simply a homothety with coefficient $\sqrt N$. The
polytope  $\De _{\sG _{ K_N }}$ is the convex hull of $\p ^* (\sum
_{i\in J} v_i)$ for all subsets $J\subset \{ 1,2, \dots ,N\}$ and
the vectors $\sum _{i\in J} v_i \in C_0 (K_N ,\R )$ are exactly the
vertices of the standard cube in $C_0 (K_N ,\R )$ defined by
inequalities $0\le a_i \le 1$,  $i=1, \ldots , N$. Therefore up to
homothety  $\De _{\sG _{ K_N }}$ is the orthogonal projection of the
standard cube of dimension $N$ parallel to its main diagonal. In
particular for $N=4$ the graph $K_4$ is  planar and
$\widehat{K_4}=K_4$, so the graphic and the  cographic systems
coincide, $\dim F _{ K_4 }=\dim G _{ K_4 }=3$ and the polyhedron
$\De _{\sF _{ K_4 }}= \De _{\sG _{ K_4 }}$ is  the rhombic
dodecahedron.
\end{ex}

\section{Geometry of the Bixby-Seymour system.}

Here we discuss geometry of the Bixby-Seymour unimodular system
(Example \ref{ex_Seymour}).
\begin{prop}\label{prop_Seymour}
\begin{enumerate}
  \item  The Bixby-Seymour system $\Psi $ is self-dual.
  \item The Bixby-Seymour system is neither graphic nor cographic.
  \item Complexity of the Bixby-Seymour system is $ c( \Psi )= 162$.
  \item The polytope $\De _{\Psi }$ is the convex hull of two
  regular simplices in Euclidean space of dimension 5 that are
  centrally symmetrical to each other. The square of the
  position vector of each vertex is 10.
  \item The lattice $L_{\Psi }$ is generated by  midpoints  of edges of
  these simplices.
  \item \label{points_on_Seymour_polytope} $\De _{\Psi } \cap L_{\Psi } $
  consists of
\begin{itemize}
  \item  12 vertices of $\De _{\Psi }$ ( vertices of the two
  simplices), the squares of the position vectors are $10$;
  \item 30  midpoints of edges of
  these simplices,   the squares of the position vectors are
  $4$;
  \item  30 midpoints of segments connecting vertices of
  different simplices,   the squares of the position vectors
  are   $6$.
\end{itemize}
\end{enumerate}
\end{prop}

The first statement is already proved in Corollary
\ref{cor_Seymour_autodual}. The third may be computed directly:
$c(\Psi ) = \det ( R^T R) $, where $R$ is the matrix
(\ref{Seymour}). The second statement can be deduced from statements
4-6 as follows. Suppose that $\Psi $ is graphic or cographic,
corresponding to a certain  graph $\Ga $. Then according to
Proposition \ref{stabilization_of_graph} the  graph  $\Ga $ has no
bridges and loops and so $\Psi $ is both graphic and cographic.
Therefore  the  graph $\Ga $ has 10 edges and 6 vertices. Each
vertex of  the  graph $\Ga $  define the cut (\ref{cut_defined
by_vertex}) which is the element of the cographic lattice and lies
on  $\De _{\Psi }$. The square of this element is equal  to the
valency of the vertex. Item \ref{points_on_Seymour_polytope} states
that such vectors may have square $4$, $6$ or $10$, and therefore
valence of each vertex is at least $4$. Thus we come to
contradiction because a graph having 6 vertices each of valency at
least 4 should have at least $(4\cdot 6) / 2 =12$ edges.

The mapping $\Phi _{\Psi }$ (\ref{def_Phi}) is defined by the matrix
(\ref{Seymour}) or by the matrix (\ref{preSeymour}) --- they differ
only in the choice of a basis in the 5-dimensional space $U$. We
shall use the matrix (\ref{preSeymour}). Note that rows of
(\ref{preSeymour}) correspond to two-element subsets of the set $\{
1,2,3,4,5 \}$, i.e. to the edges of a complete graph $K_5$, whose
vertices are labeled by this set. Let's attach variable $x_i$ to the
vertex $i$ and variable $y_{ij}$ to the edge $ij$, then the image
$\Phi _{\Psi } (U)= W _{\Psi }$ (see (\ref{def_V_Om})) is the set of
all vectors $Y = ( \ldots , y_{ij} , \ldots )$ such that $
y_{ij}=x_i+x_j$ for certain vector $X=(x_1,x_2,x_3,x_4,x_5)\in U$.
So $W _{\Psi }^{\bot }$ is the set of vectors $Z= ( \ldots ,z_{ij},
\ldots )$ such that $\sum z_{ij}(x_i+x_j)=0$ for all vectors $X\in
U$. For vectors $Z \in \De _{\Psi ^{\bot }} \cap L_{\Psi ^{\bot }} $
the values of $z_{ij}$ may be only $1$, $-1$  or $0$. Since $\Psi $
is self-dual we get the following characterization of vectors from
$\De _{\Psi  } \cap L_{\Psi   } $: these are vectors $Z=  ( \ldots ,
z_{ij} ,\ldots )$,  such that all $z_{ij}\in \{ 0, 1 ,-1\}$ and
\begin{equation}\label{z_ij}
\sum z_{ij}(x_i+x_j)=0
\end{equation}
 for all vectors $(x_1,x_2,x_3,x_4,x_5)\in
U$. For each relation (\ref{z_ij}) let us denote by $E'$ the set of
edges $(ij)$ of $K_5$ such that $z_{ij}=1$ and denote by $E''$ the
set of edges $(pq)$ of $K_5$ such that $z_{pq}=-1$. Then the
relation (\ref{z_ij}) may be rewritten as follows:
\begin{equation}\label{E'_E''}
\sum  _{(ij)\in E'}(x_i+x_j)= \sum  _{(pq)\in E''}(x_p+x_q).
\end{equation}
Denote by $\Ga '$ and   $\Ga ''$ the subgraphs of $K_5$ such that
 $E(\Ga ')=E'$ and  $E(\Ga '')=E''$. The relation (\ref{E'_E''})
holds for all $x_i$ if and only if each variable  $x_i$ occurs the
same number of times on the left and  on the right sides of the
equation. Geometrically it means that each vertex of $K_5$ has equal
valency in $\Ga '$ and   $\Ga ''$. So we see that the lattice points
of the polytope   $\De _{\Psi  } $ are in one to one correspondence
with pairs of subgraphs $\Ga '$ and   $\Ga ''$ of $K_5$ such that
 $\Ga '$ and   $\Ga ''$ have no common edges and each
 vertex of $K_5$ has equal valency
in $\Ga '$ and in  $\Ga ''$. Now it is quite easy to list such pairs
of subgraphs,  especially considering that  valences of vertices in
each subgraph cannot be greater than  $2$. Only three such
configurations are possible, they are shown in Figures 1-3. The
edges of the graph $\Ga '$ are shown there by solid lines, and the
edges of the  graph $\Ga ''$ are shown by dotted lines.

\begin{tabular}{ccc}
\begin{picture}(100,100) 
\put(30,10){\circle*{5}} \put(70,10){\circle*{5}}
\put(10,50){\circle*{5}} \put(90,50){\circle*{5}}
\put(50,90){\circle*{5}}
\multiput(30,10)(3,2){20} {\circle*{2}}
\multiput(10,50)(3,-2){20} {\circle*{2}}
\put(30,10){\line(-1, 2){20}}
\put(70,10){\line(1, 2){20}}
\end{picture}
&
\begin{picture}(100,100)  
\put(30,10){\circle*{5}} \put(70,10){\circle*{5}}
\put(10,50){\circle*{5}} \put(90,50){\circle*{5}}
\put(50,90){\circle*{5}}
\put(10,50){\line( 1, 1){40}}
\put(90,50){\line( -1, 1){40}}
\put(30,10){\line( 1, 0){40}}
\multiput(30,10)(1,4){20} {\circle*{2}}
\multiput(70,10)(-1,4){20} {\circle*{2}}
\multiput(10,50)( 4,0){20} {\circle*{2}}
\end{picture}
&
\begin{picture}(100,100)    
\put(30,10){\circle*{5}} \put(70,10){\circle*{5}}
\put(10,50){\circle*{5}} \put(90,50){\circle*{5}}
\put(50,90){\circle*{5}}
\put(30,10){\line(-1, 2){20}}
\put(70,10){\line(1, 2){20}}
\put(10,50){\line( 1, 1){40}}
\put(90,50){\line( -1, 1){40}}
\put(30,10){\line( 1, 0){40}}
\multiput(30,10)(1,4){20} {\circle*{2}}
\multiput(30,10)(3,2){20} {\circle*{2}}
\multiput(10,50)(3,-2){20} {\circle*{2}}
\multiput(70,10)(-1,4){20} {\circle*{2}}
\multiput(10,50)( 4,0){20} {\circle*{2}}
\end{picture}
\\
Figure 1 &  Figure 2  &  Figure 3\\
\end{tabular}

Each such picture corresponds to a point  $Z= ( \ldots , z_{ij} ,
\ldots ) \in \De _{\Psi } \cap L_{\Psi } $, where the solid edge
$(ij)$ correspond to the coordinate  $z_{ij}=1$, the dotted edge
$(pq)$ correspond to the coordinate  $z_{pq}=-1$,  and the missing
edges correspond to zero values of   $z_{ij}$. So the squares of the
vectors corresponding to figures 1, 2 and 3 are respectively $4$,
$6$ and $10$. The permutation group $S_5$ acts on the pictures of
each type and the stabilizers of the figures 1, 2 and 3 are
subgroups of  $S_5$ of order, respectively $4$, $4$ and $10$.
Therefore the orbits of the configurations from figures
 1, 2 and 3 have respectively $30$,
$30$ and $12$ elements --- these are the numbers of corresponding
lattice points of  $\De _{\Psi }$. Next we can calculate the scalar
products of vectors corresponding to different configurations from
figure 3 and verify that they form the vertices of two regular
simplices that are centrally symmetrical to each other. The vertices
of each of the simplices are the orbits of the alternating group
$A_5\subset S_5$. To complete the proof of item
\ref{points_on_Seymour_polytope}, we need to represent the vectors
corresponding to figures 1 and 2  as  half-sums of the vectors
corresponding to figure 3. These representations are presented in
figures 4 and 5.

$$
\begin{picture}(100,50)(0,50) 
\put(30,10){\circle*{5}} \put(70,10){\circle*{5}}
\put(10,50){\circle*{5}} \put(90,50){\circle*{5}}
\put(50,90){\circle*{5}}
\multiput(30,10)(3,2){20} {\circle*{2}}
\multiput(10,50)(3,-2){20} {\circle*{2}}
\put(30,10){\line(-1, 2){20}}
\put(70,10){\line(1, 2){20}}
\end{picture}
= \frac 12 \left (
\begin{picture}(100,50)(0,50)    
\put(30,10){\circle*{5}} \put(70,10){\circle*{5}}
\put(10,50){\circle*{5}} \put(90,50){\circle*{5}}
\put(50,90){\circle*{5}}
\put(30,10){\line(-1, 2){20}}
\put(70,10){\line(1, 2){20}}
\put(10,50){\line( 1, 1){40}}
\put(90,50){\line( -1, 1){40}}
\put(30,10){\line( 1, 0){40}}
\multiput(30,10)(1,4){20} {\circle*{2}}
\multiput(30,10)(3,2){20} {\circle*{2}}
\multiput(10,50)(3,-2){20} {\circle*{2}}
\multiput(70,10)(-1,4){20} {\circle*{2}}
\multiput(10,50)( 4,0){20} {\circle*{2}}
\end{picture}
+
\begin{picture}(100,50)(0,50)    
\put(30,10){\circle*{5}} \put(70,10){\circle*{5}}
\put(10,50){\circle*{5}} \put(90,50){\circle*{5}}
\put(50,90){\circle*{5}}
\put(30,10){\line(-1, 2){20}}
\put(70,10){\line(1, 2){20}}
\multiput(10,50)(4,4){10} {\circle*{2}}
\multiput(50,90)(4,-4){10} {\circle*{2}}
\multiput(30,10)( 4,0){10} {\circle*{2}}
\put(30,10){\line( 1, 4){20}}
\multiput(30,10)(3,2){20} {\circle*{2}}
\multiput(10,50)(3,-2){20} {\circle*{2}}
\put(70,10){\line(-1, 4){20}}
\put(10,50){\line( 1, 0){80}}
\end{picture}
\right )
$$

\begin{center}Figure 4
\end{center}

$$
\begin{picture}(100,50)(0,50)  
\put(30,10){\circle*{5}} \put(70,10){\circle*{5}}
\put(10,50){\circle*{5}} \put(90,50){\circle*{5}}
\put(50,90){\circle*{5}}
\put(10,50){\line( 1, 1){40}}
\put(90,50){\line( -1, 1){40}}
\put(30,10){\line( 1, 0){40}}
\multiput(30,10)(1,4){20} {\circle*{2}}
\multiput(70,10)(-1,4){20} {\circle*{2}}
\multiput(10,50)( 4,0){20} {\circle*{2}}
\end{picture}
= \frac 12 \left (
\begin{picture}(100,50)(0,50)    
\put(30,10){\circle*{5}} \put(70,10){\circle*{5}}
\put(10,50){\circle*{5}} \put(90,50){\circle*{5}}
\put(50,90){\circle*{5}}
\put(30,10){\line(-1, 2){20}}
\put(70,10){\line(1, 2){20}}
\put(10,50){\line( 1, 1){40}}
\put(90,50){\line( -1, 1){40}}
\put(30,10){\line( 1, 0){40}}
\multiput(30,10)(1,4){20} {\circle*{2}}
\multiput(30,10)(3,2){20} {\circle*{2}}
\multiput(10,50)(3,-2){20} {\circle*{2}}
\multiput(70,10)(-1,4){20} {\circle*{2}}
\multiput(10,50)( 4,0){20} {\circle*{2}}
\end{picture}
+
\begin{picture}(100,50)(0,50)    
\put(30,10){\circle*{5}} \put(70,10){\circle*{5}}
\put(10,50){\circle*{5}} \put(90,50){\circle*{5}}
\put(50,90){\circle*{5}}
\multiput(70,10)(2,4){10} {\circle*{2}}
\multiput(30,10)(-2,4){10} {\circle*{2}}
\put(10,50){\line( 1, 1){40}}
\put(90,50){\line( -1, 1){40}}
\put(30,10){\line( 1, 0){40}}
\multiput(30,10)(1,4){20} {\circle*{2}}
\multiput(70,10)(-1,4){20} {\circle*{2}}
\put(30,10){\line( 3, 2){60}}
\put(70,10){\line( -3, 2){60}}
\multiput(10,50)( 4,0){20} {\circle*{2}}
\end{picture}
\right )
$$
\begin{center}Figure 5
\end{center}

Thus the 60 midpoints are not the vertices of $\De _{\Psi }$, so the
only vertices of $\De _{\Psi }$ are the 12 vertices of the two
simplices, which proves item 4). For item 5) it is sufficient to
notice that the vectors with square 4 are the shortest vectors in
the lattice $L_{\Psi }$ and therefore any five such linearly
independent vectors form a basis of the lattice.  If a vector  $Z= (
\ldots ,z_{ij}, \ldots )\in L_{\Psi }$ has square $1$, $2$ or $3$
then all the entries $z_{ij}$ should be $0$, $1$ or $-1$ and
therefore such $ Z\in L_{\Psi } \cap \De _{\Psi } $, but we have
seen that squares of vectors from $ L_{\Psi } \cap \De _{\Psi } $
may be only  $4$, $6$ or $10$. $\square $

The automorphisms of the Bixby-Seymour
  system are the permutations of the vertices of the simplex and
the central symmetry that interchanges the two simplices.

\begin{cor}\label{cor_AUT}
The automorphism group of the Bixby-Seymour
  system is ${\Z _2 \times S_6}$.
\end{cor}

Note that our approach also allows us to describe 20  facets (see
Proposition \ref{faces_polytope}) of the polytope $\De _{\Psi } $.
10 facets are defined by equations $z_{ij}=1$ and so the lattice
points of corresponding facets are given by all the configurations
from Figures 1-3, containing the solid edge $(ij)$. Similarly 10
other facets
 are defined by equations $z_{ij}=-1$ and the lattice
points of these facets correspond to all configurations from Figures
1-3, containing the dotted edge $(ij)$. It is not hard to verify
that each facet has exactly 6 vertices: three from one simplex and
the remaining three are vertices of the second simplex, not opposite
to the first three. So the facets of $\De _{\Psi } $  are in
one-to-one correspondence with the triplets of vertices of one of
the simplices.

\end{document}